\documentclass[12pt,a4paper]{amsart}
\usepackage[top=2.5cm,bottom=2.5cm,left=2.5cm,right=2.5cm]{geometry}
\usepackage{cases}
\usepackage{graphicx}
\usepackage{amssymb}
\usepackage{mathrsfs}
\usepackage{amsthm}
\usepackage{amsmath}
\usepackage{ulem}
\usepackage{xcolor}
\newtheorem{problem}{Problem}
\newtheorem{theorem}{Theorem}[section]
\newtheorem{corollary}[theorem]{Corollary}
\newtheorem{lemma}[theorem]{Lemma}

\newtheorem{definition}[theorem]{Definition}
\newtheorem{remark}{Remark}

\numberwithin{equation}{section}
\newtheorem{example}{Example}

\newcommand{\beqm}{\begin{eqnarray*}}
	\newcommand{\eqm}{\end{eqnarray*}}

\newcommand{\C}{\mathbb{C}}

\newcommand{\D}{\mathbb{D}}

\begin{document}
\title[Schwarz-Pick estimates of nonvanishing holomorphic functions]{Schwarz-Pick estimates of nonvanishing holomorphic functions and their application}
\author{Junquan Pan}
\address{Department of Mathematics, Huzhou Normal University, Huzhou, Zhejiang 313000 China}
\email{pjq13626602189@163.com}

\author{Xiaomin Tang}
\address{Department of Mathematics, Huzhou Normal University, Huzhou, Zhejiang 313000 China}
\email{txm@huznu.edu.cn}

\author{Wei Chen*}
\address{Department of Mathematics, Huzhou Normal University, Huzhou, Zhejiang 313000 China}
\email{chenwei0721@huznu.edu.cn}
\thanks {*Corresponding author}
\begin{abstract}
 The main purpose of this paper is to study Schwarz-Pick estimates for holomorphic functions avoiding the origin. We first prove higher order derivative estimates for holomorphic functions from the unit disk $\mathbb{D}$ to the punctured disk. Then we extend this result to the unit ball $\mathbb{B}_n$ in $\mathbb{C}^n$, and establish a Schwarz-Pick type inequality for Fr\'{e}chet derivatives of holomorphic functions. Finally, we apply our new Schwarz-Pick estimates to obtain a $\L$ojasiewicz-type inequality  and give the explicit expression of the $\L$ojasiewicz exponent.
\end{abstract}

\maketitle
\section{Introduction}
Let $\mathbb{D}$ be the unit disk in the complex plane $\mathbb{C}$, and let $\mathbb{B}_n$ be the unit ball in the $n$-dimensional complex
Euclidean space $\mathbb{C}^n$. For any $z,w\in\mathbb{C}^n$, we denote $\langle z,w\rangle=z_1\overline{w}_1+\cdots+z_n\overline{w}_n$
and $|z|=\sqrt{|z_1|^2+\cdots+|z_n|^2}$.
A multi-index $\alpha=\left(\alpha_{1},\cdots,\alpha_{n}\right)$ consists of $n$ nonnegative integers $\alpha_{j}$ for $1\leq j\leq n$.
The degree and factorial of $\alpha$ are given by
$\left|\alpha\right|=\sum_{j=1}^{n}\alpha_{j}$ and  $\alpha!=\alpha_{1}!\cdots \alpha_{n}!$,
respectively.\par
 The classical Schwarz lemma states that holomorphic self-maps of $\mathbb{D}$ have the contraction property. It is well known that Schwarz lemma has become a crucial theme in many branches of mathematical research for more than a hundred years to date. Here we refer the reader to \cite{HXJ1,HXJ2, LP, LT, RB, TLZ, Zh25} for more on these matters. And then it is natural to consider estimates for derivatives of holomorphic functions. There is the classical Schwarz-Pick lemma as follows.
\begin{theorem}[\cite{Pick}]
Let $f:\mathbb{D}\longrightarrow\mathbb{D}$ be a holomorphic function. Then for any $z\in \mathbb{D}$,
$$|f'(z)|\leq \frac{{1-|f(z)|^2}}{1-|z|^2}.$$
\end{theorem}

A further profound generalization is Ahlfors-Schwarz lemma that extends Schwarz-Pick lemma to Riemannian surface from the viewpoint of differential geometry.
Ahlfors \cite{Ahlfors} proved that any holomorphic mapping from the unit disk into the Riemann surface with the Hermitian metric of Gauss curvature at most $-4$ is distance-decreasing with respect to Poincar$\acute{\text{e}}$ metric.
More recently, attention has been paid to Schwarz-Pick lemma in different domains. See \cite{BKR,CH,FK,GKMR, LC,PR} for more on
these matters. Naturally, the classical Schwarz-Pick estimate can be generalized to derivatives of arbitrary order for holomorphic functions on $\mathbb{D}$.
See the relevant works \cite{AW2003,GZ2004,MSZ2002}. Specially, Dai and Pan obtained the following more explicit Schwarz-Pick estimate in \cite{DP2008}.
\begin{theorem}[\cite{DP2008}]
Let $f: \mathbb{D}\longrightarrow \mathbb{D}$ be a holomorphic function. Then for any $z\in\mathbb{D}$ and $k\geq 1$,
$$
\left|f^{(k)}(z)\right|\leq\frac{k!(1-\left|f(z)\right|^2)}{(1-\left|z\right|^2)^k}(1+|z|)^{k-1}.
$$
\end{theorem}
A natural question arises: is it possible to provide Schwarz-Pick estimates for holomorphic functions when image domains are  multiply connected punctured domains? In this paper, we first study Schwarz-Pick estimates for holomorphic functions from $\mathbb{D}$ into punctured domains. One of our main results is stated as follows.
\begin{theorem}\label{cw2025}
	Let $f:\mathbb{D}\rightarrow\mathbb{D}\setminus\left\{0\right\}$ be a holomorphic function. Then for each $z\in\D$ and $m\geq 1$,
 \begin{equation}\label{20251224}\begin{split}
	\left|f^{(m)}(z)\right|\leq \frac{|f(z)|}{\left(1-\left|z\right|\right)^{m}}
\sum\limits_{k=1}^{m}\frac{m!}{k!} \! \begin{pmatrix}
m-1 \\
k-1
\end{pmatrix}\left(\frac{2\log\frac{1}{|f(z)|}}{1+\left|z\right|}\right)^k.
	\end{split}\end{equation}
\end{theorem}
Let $E$ and $F$ be two complex Banach spaces with norms $\|\cdot\|_{E}$ and $\|\cdot\|_{F}$, respectively. Let $U\subseteq E$ be an open subset. A continuous mapping $f:U\longrightarrow F$ is said to be holomorphic in $U$ if for any $z\in U$, there exists a bounded linear operator  $A_z:E\longrightarrow F$ such that
	$$
	\lim_{\|h\|_E\rightarrow0}\frac{\|f\left(z+h\right)-f\left(z\right)-A_z\cdot h\|_{F}}{\|h\|_{E}}=0.
 	$$
  The linear operator $A_z$ is called the Fr$\acute{\text{e}}$chet derivative of $f$ at $z$ and is written as $Df\left(z\right)$. Also, $Df\left(z\right)\cdot h$ is called the Fr$\acute{\text{e}}$chet derivative of $f$ at $z$ in the direction $h\in E$.
  Furthermore, if $f: U\longrightarrow F$ is a holomorphic mapping, then
  $$
  f\left(w\right)=\sum_{k=0}^{\infty}\frac{1}{k!}D^{k}f\left(z\right)\cdot\left(w-z\right)^k
  $$
  for all $w$ in some neighborhood of $z\in U$, where $D^kf(z)$ is the $k$th Fr$\acute{\text{e}}$chet derivative of $f$ at $z$, and for $k\geq1$,
  $$
  D^{k}f\left(z\right)\cdot\left(w-z\right)^k=D^kf(z)\cdot(w-z,\cdots,w-z).
  $$
  Moreover, $D^kf(z)$ is a bounded symmetric $k$-multilinear mapping from $\prod\limits_{j=1}^{k}E$ into $F$.
  In particular, if $E=\mathbb{C}^n$ and $F=\mathbb{C}^m$, then
 $$
 D^{k}f(z)\cdot \beta^k=\sum\limits_{|\alpha|=k}\frac{k!}{\alpha!}\frac{\partial^kf(z)}{\partial z_{1}^{\alpha_{1}}\cdots\partial z_{n}^{\alpha_{n}}}\beta^{\alpha},\ z\in U \subset\mathbb{C}^{n},\ \beta\in\mathbb{C}^{n},
 $$
 where $\alpha=(\alpha_1,\cdots,\alpha_n)$ is a multi-index.
 In this paper, the second purpose is to give the following Schwarz-Pick estimates for higher order Fr$\acute{\text{e}}$chet derivatives of nonvanishing holomorphic functions on $\mathbb{B}_n$.
\begin{theorem}\label{pjq2026-3}
 Let $f:\mathbb{B}_n\rightarrow\mathbb{D}\setminus\left\{0\right\}$ be a holomorphic function. Then for any $z\in\mathbb{B}_n$ and $\beta\in\partial \mathbb{B}_n$,
	\begin{equation}\label{cw2032}\begin{split}
		\left|D^mf(z)\cdot\beta^m\right|\leq |f(z)|\left(\frac{\sqrt{A}+\left|\langle z,\beta\rangle\right|}{1-|z|^2}\right)^m
    \sum\limits_{k=1}^{m}\frac{m!}{k!} \! \begin{pmatrix}
    m-1 \\
    k-1
    \end{pmatrix}\left(\frac{2\log\frac{1}{|f(z)|}}{\sqrt{A}+\left|\langle z,\beta\rangle\right|}\right)^k,
		\end{split}\end{equation}
where $A=1-\left|z\right|^2+\left|\langle z,\beta\rangle\right|^2$ and $m\geq1$.
\end{theorem}

In 1959, $\L$ojasiewicz \cite{LS} proved that if $f:G\rightarrow\mathbb{R}$ is a real analytic function, where $G\subset\mathbb{R}^n$ is an open set, $f^{-1}(0)\neq\emptyset$ and $E\subset G$ is a compact set, then there exists $C,\alpha>0$ such that for any $x\in E$,
$$
|f(x)|\geq C\left(dist(x;f^{-1}(0))\right)^\alpha.
$$
This deep result later was called $\L$ojasiewicz inequality, which shows the profound relationship between real analytic functions and their zero sets. $\L$ojasiewicz inequality has been studied in real algebraic geometry and partial differential equations. Here we refer the reader to \cite{MA2022,PR2018}. The \L ojasiewicz exponent $\alpha$ is critical in that it describes the decay rate of $f$ tending to 0 as $x$ approaches the zero set. The classical \L ojasiewicz inequality ensures the existence of such an exponent, but no explicit expression for this exponent $\alpha$ has been obtained. To work out this explicit form of the $\L$ojasiewicz exponent is then an interesting problem. It is natural to consider the $\L$ojasiewicz-type inequality in the case of several complex variables. To proceed further, we first introduce some notation and definitions.

Let $\Omega$ be a domain in $\mathbb{C}^n\ (n\geq2)$. For a complex function $f$ on $\Omega$ and $p\in\Omega$, if there exists a multi-index $\alpha=(\alpha_1,\cdots,\alpha_{n})$ with $|\alpha|=k$ satisfying
$$ \frac{\partial^{|\alpha|} f}{\partial z^{\alpha_{1}}_1\cdots\partial z_{n}^{\alpha_n}}(p) \neq 0, \ \frac{\partial^{|\beta|} f}{\partial z^{\beta_{1}}_1\cdots\partial z_{n}^{\beta_n}}(p)=0$$
for any multi-index $\beta=(\beta_1,\cdots,\beta_n)$ with $|\beta|<k$, then we call $k$ the vanishing order of $f$ at $p$. Let $Z_f=\{z\in\Omega:f(z)=0\}$ be the zero set of $f$ in $\Omega$ and $\mathbf{0}=(0,\cdots,0)$. In this paper, the third purpose is to apply the Schwarz-Pick estimate to obtain a $\L$ojasiewicz-type inequality.

\begin{theorem}\label{wst2030}
	Let $\Omega\subset\mathbb{C}^n\ (n\geq2)$ be a domain with $q\in\Omega$. Suppose that $f$ is a holomorphic function on $\Omega$ and $f(q)=0$. If $k$ is the vanishing order of $f$ at $q$, then there exists a neighborhood $U$ of $q$ such that for any $0<\varepsilon<1$ and $z\in U$,
	$$
	|f(z)|\geq C\operatorname{dist}\left(z, Z_{f}\right)^{\frac{k}{1-\varepsilon}},
	$$
	where the constant $C$ depends on $f$ and $\varepsilon$.
\end{theorem}	

As an application of the above Schwarz-Pick estimate, we characterize the behavior of holomorphic functions near their zero sets. In contrast to the general real analytic case where the exponent remains implicit, we give an explicit expression for the $\L$ojasiewicz exponent. However, it should be noted that this expression fails to hold for real analytic functions, as shown by the following counterexample.
\begin{example}
For $m\geq1$, let $f(x,y)=x^2+y^{2m}$, where $(x, y)\in \mathbb{R}^2$. Then $Z_f=f^{-1}(0)=\{(0, 0)\}$ and $\text{dist}((x,y), Z_f)=\sqrt{x^2+y^2}$. Also, the vanishing order of $f$ is $2$. Suppose
$$\left|f(x,y)\right|\geq C(\sqrt{x^2+y^2})^{\frac{2}{1-\varepsilon}}.$$
Taking $x=0$, this leads to a contradiction for $0<\varepsilon<1-\frac{1}{m}$.
\end{example}
The rest of the paper is organized as follows. In Section 2, we provide several higher order derivative estimates for nonvanishing holomorphic functions on $\mathbb{D}$. In Section 3, we give Schwarz-Pick estimates for Fr$\acute{\text{e}}$chet derivatives of holomorphic functions from $\mathbb{B}_n$ into the punctured disk. In Section 4,
we are devoted to establishing the $\L$ojasiewicz-type inequality in several complex variables setting.

\section{Schwarz-Pick Estimates on the unit disk}
In this section, we first introduce some lemmas, which not only are used in the subsequent sections but also have their interest.
And then we prove Theorem 1.3.
\begin{lemma}[Fa$\grave{\text{a}}$ di Bruno's formula, \cite{Ch02, Co74}]\label{2.1} If $g$ and $f$ are complex functions with a sufficient number of derivatives, then
	$$\frac{d^m}{dz^m}f\left(g(z)\right)=\sum_{k=1}^mf^{(k)}\left(g(z)\right) B_{m, k}\left(g'(z),g''(z),\cdots,g^{(m-k+1)}(z)\right),$$
	where the Bell polynomial $B_{m,k}$ is defined by
	$$B_{m, k}(z_1,z_2,\cdots,z_{m-k+1})=\sum_{\nu_1,\cdots,\nu_{m-k+1}}\frac{m!}{\nu_1!\nu_2!\cdots \nu_{m-k+1}!}\prod_{j=1}^{m-k+1}\left(\frac{z_j}{j!}\right)^{\nu_j}.$$
The sum is over all sequences of nonnegative integers $\left\{\nu_{1}, \nu_{2}, \cdots, \nu_{m-k+1}\right\}$ such that the following two conditions are satisfied:
 $\nu_{1}+2\nu_{2}+\cdots+(m-k+1)\nu_{m-k+1}=m$ and $\nu_{1}+\nu_{2}+\cdots+\nu_{m-k+1}=k$. In particular, it holds that
 \begin{equation}\label{pan5}B_{m,k}\left(1!,2!,\cdots,(m-k+1)!\right)=
\frac{m!}{k!} \! \begin{pmatrix}
m-1 \\
k-1
\end{pmatrix}\end{equation}
and for any $ a,b \in \mathbb{C} $,
\begin{equation}\label{pan3}
B_{m,k}(abz_1,\, ab^2z_2,\, \cdots,\, ab^{m-k+1}z_{m-k+1})
= a^k b^m\, B_{m ,k}(z_1,\, z_2,\, \cdots,\, z_{m-k+1}).
\end{equation}
		
\end{lemma}
\begin{lemma}[Riesz-Herglotz theorem, \cite{AL1979}]\label{RH}
A function $f$ on $\mathbb{D}$ is holomorphic and has non-negative real part if and only if there exists a finite positive Borel measure $\mu$ on $[0,2\pi)$ such that for any $z\in\mathbb{D}$,
	\begin{align*}
	f(z) = i \text{Im} f(0) + \frac{1}{2\pi} \int_{[0, 2\pi)} \frac{e^{it} + z}{e^{it}-z} \, d\mu(t).
	\end{align*}
\end{lemma}
The following lemma is well-known and the starting point for our study. We refer the reader to \cite{Graham}. Here, we give the proof for completeness.
\begin{lemma}\label{lemma1}
	Let	$f:\mathbb{D}\rightarrow\mathbb{H}^{+}$ be a holomorphic function with the power series $$f(z)=\displaystyle \sum_{k=0}^{\infty}a_{k}z^{k},$$ where $\mathbb{H}^{+}=\left\{z\in \mathbb{C}:\mathrm{Re}{z}>0\right\}$. Then for any $k\geq1$,
	\begin{align*}
	|a_k| \leq 2 \mathrm{Re} a_0.
	\end{align*}
\end{lemma}
\begin{proof}
	Note that
	$$	
	\frac{e^{it} + z}{e^{it}-z}= 1 + 2\sum_{k = 1}^{\infty}e^{-ikt}z^{k}
	$$
	for any $z \in \mathbb{D}$. By Lemma \ref{RH}, we have
	\begin{align*}
		f(z)
		&= i \text{Im} f(0) + \frac{1}{2\pi} \int_{[0, 2\pi)} \frac{e^{it} +z}{e^{it} - z} d\mu(t) \\
		&= i \text{Im} f(0) + \frac{1}{2\pi} \int_{[0,2\pi)} \left(1 + 	2\sum_{k=1}^{\infty} e^{-ikt}z^k\right) d\mu(t) \\
		&= i \text{Im} f(0)+ \frac{1}{2\pi}  \mu([0,2\pi)) + \frac{2}{2\pi} \int_{[0,2\pi)} \sum_{k=1}^{\infty} e^{-ikt} z^k d\mu(t) \\
		&= i \text{Im} f(0)+ \frac{1}{2\pi} \mu([0,2\pi)) + \frac{2}{2\pi} \sum_{k=1}^{\infty} \left( \int_{[0,2\pi)} e^{-ikt} d\mu(t) \right) z^k
	\end{align*}
	for any $z\in\mathbb{D}$. Hence, we obtain
	$$	
    a_0=i \text{Im} f(0) + \frac{1}{2\pi} \mu([0,2\pi))
    $$
    and
	$$	
    a_k=\frac{2}{2\pi} \int_{[0,2\pi)} e^{-i k t} d\mu(t),\ \ k\geq1.
    $$
	Taking the modulus on both sides of $a_k$ implies
	\begin{align*}
		|a_{k}|
		&=\left|\frac{2}{2\pi}\int_{[0,2\pi)}e^{-ikt}d\mu(t)\right|\leq\frac{2}{2\pi}\int_{[0,2\pi)}|e^{-ikt}|d\mu(t)\leq2\cdot\frac{1}{2\pi}\mu([0,2\pi))=2\text{Re}a_{0}.
	\end{align*}
 The proof is complete.
\end{proof}

\begin{lemma}\label{im}
Let $\phi(z)=\sum\limits_{k=0}^{\infty}a_k\left(\frac{z-\xi}{1-\overline{\xi}z}\right)^k$ be a holomorphic function on $\mathbb{D}$ for given $\xi\in\mathbb{D}$. If there exists $r_0>0$ such that $|a_k|\leq r_0$ for any $k\geq1$, then for any $m\geq 1$,
$$|\phi^{(m)}(\xi)|\leq \frac{m!r_0(1+|\xi|)^{m-1}}{(1-|\xi|^2)^m}.$$
\end{lemma}

\begin{proof}
    Note that $(\frac{z-\xi}{1-\overline{\xi}z})^k$ has a non-zero $m$th derivative at $z=\xi$ only when $k \leq m$. Then
    $$
    \phi^{(m)}(\xi)=\left(\sum\limits_{k=0}^{\infty}a_k\left(\frac{z-\xi}{1-\overline{\xi}z}\right)^k\right)^{(m)}\bigg|_{z=\xi}=\left(\sum\limits_{k=1}^{m}a_k\left(\frac{z-\xi}{1-\overline{\xi}z}\right)^k\right)^{(m)}\bigg|_{z=\xi}.
    $$
    Let us consider the derivative of $\left(\frac{z-\xi}{1-\overline{\xi}z}\right)^k$. Set
\begin{center}
    $f(z)=(z-\xi)^{k}$\ and\ $g(z)=(1-\overline{\xi}z)^{-k}$.
\end{center}
    Then by Leibniz's rule, the $m$th derivative of $f(z)g(z)$ is given by
    $$
    \frac{d^{m}}{d z^{m}}\left(f(z)g(z)\right)\bigg|_{z=\xi} =\sum_{s=0}^{m}\binom{m}{s}f^{(s)}(\xi)g^{(m-s)}(\xi).
    $$
   It is easy to check that $f^{(s)}(\xi) = s!$ if $s=k$, and $f^{(s)}(\xi) = 0$ if $s \neq k$. This implies that, for $1\leq k \leq m$,
 \begin{equation}\label{pjq2025}
    \frac{d^{m}}{d z^{m}}\left(f(z)g(z)\right)\bigg|_{z=\xi} =\binom{m}{k}f^{(k)}(\xi)g^{(m-k)}(\xi) = \binom{m}{k}k!g^{(m-k)}(\xi).
 \end{equation}
    Now we need to compute $g^{(m-k)}(\xi)$. A straightforward computation for $l$th derivative of $g(z)$ yields
    $$g^{(l)}(z)=(-k)(-k-1)\cdots(-k-l+1)(-\overline{\xi})^l (1-\overline{\xi}z)^{-k-l} = \frac{k(k+1)\cdots(k+l-1)}{(1-\overline{\xi}z)^{k+l}}(\overline{\xi})^l.$$
    Setting $l=m-k$, we have
    $$g^{(m-k)}(\xi) = \frac{k(k+1)\cdots(m-1)}{(1-\overline{\xi}\xi)^{k+(m-k)}}(\overline{\xi})^{m-k} = \frac{(m-1)!/(k-1)!}{(1-|\xi|^2)^m}(\overline{\xi})^{m-k}.$$
    This, together with (\ref{pjq2025}), shows
 \begin{align*}
    \frac{d^{m}}{d z^{m}}\left(\frac{z-\xi}{1-\overline{\xi}z}\right)^k\bigg|_{z=\xi}&= \binom{m}{k}k! \frac{(m-1)!/(k-1)!}{(1-|\xi|^2)^m}(\overline{\xi})^{m-k}\\
     &= \frac{m!}{k!(m-k)!}k! \frac{(m-1)!}{(k-1)!} \frac{(\overline{\xi})^{m-k}}{(1-|\xi|^2)^m}\\
     &=\frac{m!}{(m-k)!} \frac{(m-1)!}{(k-1)!} \frac{(\overline{\xi})^{m-k}}{(1-|\xi|^2)^m}
 \end{align*}
    for $1\leq k \leq m$. Combining these results gives
    $$
    \phi^{(m)}(\xi)=\sum_{k=1}^{m} a_{k} \frac{m!}{(m-k)!} \frac{(m-1)!}{(k-1)!} \frac{(\overline{\xi})^{m-k}}{(1-|\xi|^2)^m}.
    $$
    It follows from $|a_k|\leq r_0$ that
  \begin{align*}
     |\phi^{(m)}(\xi)|
    &\leq \sum_{k=1}^{m}|a_{k}| \left| \frac{m!}{(m-k)!} \frac{(m-1)!}{(k-1)!} \frac{(\overline{\xi})^{m-k}}{(1-|\xi|^2)^m} \right| \\
    &\leq \sum_{k=1}^{m}r_0 \frac{m!}{(m-k)!} \frac{(m-1)!}{(k-1)!} \frac{|\xi|^{m-k}}{(1-|\xi|^2)^m} \\
    &= \frac{m!r_0}{(1-|\xi|^2)^m} \sum_{k=1}^{m} \frac{(m-1)!}{(m-k)!(k-1)!}|\xi|^{m-k}.
\end{align*}
    Notice that the binomial expansion of $(1+|\xi|)^{n-1}$ is
   $$
   (1+|\xi|)^{m-1}=\sum_{k=1}^{m} \frac{(m-1)!}{(m-k)!(k-1)!}|\xi|^{m-k}.
   $$
   Therefore, we obtain
   $$
   |\phi^{(m)}(\xi)| \leq \frac{m!r_0(1+|\xi|)^{m-1}}{(1-|\xi|^2)^m}.
   $$
 The proof is complete
 \end{proof}
    \textbf{Proof of Theorem \ref{cw2025}.} Take $g(z)=-\log f(z)$. Since $g$ maps a simply connected domain to a simply connected domain, we can take the principal branch of $g$. Hence, $g:\mathbb{D}\rightarrow \mathbb{H}^{+}$ is holomorphic, where $\mathbb{H}^+=\left\{z\in\mathbb{C}:\text{Re} z>0\right\}$.
	Thus, we have $f(z)=e^{-g(z)}$ with
    $$
    \text{Re} g(z)=\log\frac{1}{|f(z)|}>0.
    $$
    For any given $\xi\in\mathbb{D}$, let $\varphi_\xi:\mathbb{D}\longrightarrow\mathbb{D}$ be the M$\ddot{\text{o}}$bius transformation, defined by
	$$
	\varphi_\xi\left(z\right)=\frac{\xi+z}{1+\overline{\xi}z}.
	$$
	Take
	$$G(z)=g\circ \varphi_{\xi}(z)=\sum_{l=0}^{\infty}a_{l}z^{l}, \ z\in \mathbb{D}.$$
	Then $G:\mathbb{D}\rightarrow \mathbb{H}^{+}$ is a holomorphic function. From Lemma \ref{lemma1},
    \begin{equation}\label{pan1}
    |a_l|\leq2\text{Re}{a_0}=2\text{Re}{G(0)}=2\text{Re}{g(\xi)}=2\log\frac{1}{\left|f(\xi)\right|}
    \end{equation}
    for $l\geq1$.\ It follows that
	 $$f(z)=e^{-g(z)}=e^{-G(h(z))},$$
 where $h(z)=\frac{z-\xi}{1-\overline\xi z}$ is the inverse mapping of $\varphi_{\xi}$.
	By Lemma \ref{2.1}, we obtain
	\begin{align*}
   f^{(m)}(z)=&\sum_{k=1}^{m}(-1)^{k}e^{-G(h(z))}\\
   &\times B_{m,k}\left(\left(\sum_{l=0}^{\infty}a_{l}\left(h\left(z\right)\right)^{l}\right)',\left(\sum_{l=0}^{\infty}a_{l}\left(h\left(z\right)\right)^{l}\right)'',\cdots,\left(\sum_{l=0}^{\infty}a_{l}\left(h\left(z\right)\right)^{l}\right)^{(m-k+1)}\right)
    \end{align*}
	for $m\geq1$, where $B_{m,k}$ is the Bell polynomial. Note that
    \begin{center}
    $a_{0}=G(0)=g(\xi)=-\log f(\xi)$\ \  and\ \ $G(h(\xi))=G(0)$.
    \end{center}
   This implies
    \begin{equation}\label{im1}\begin{split}
   f^{(m)}(\xi)=\sum_{k=1}^{m}(-1)^{k}f(\xi)B_{m,k}\left(\left(\sum_{l=0}^{\infty}a_{l}\left(h\left(z\right)\right)^{l}\right)',\cdots,\left(\sum_{l=0}^{\infty}a_{l}\left(h\left(z\right)\right)^{l}\right)^{(m-k+1)}\right)\bigg|_{z=\xi}.
    \end{split}\end{equation}
    Moreover, by Lemma \ref{im} and (\ref{pan1}), we get
    \begin{equation}\label{pan2}
   \left|\left(\sum_{l=0}^{\infty}a_{l}\left(h\left(z\right)\right)^{l}\right)^{(k)}\bigg|_{z=\xi}\right|\leq \frac{k!\left(2\log \frac{1}{|f(\xi)|}\right)}{(1-|\xi|^2)^k}(1+|\xi|)^{k-1}.
    \end{equation}
    It follows from (\ref{pan3}), (\ref{im1}) and (\ref{pan2}) that
    \begin{equation}\label{pan4}\begin{split}
	\left|f^{(m)}(\xi)\right|&\leq |f(\xi)|\\
&\times\sum\limits_{k=1}^{m}B_{m,k}\left(\frac{1!}{1-|\xi|^2},\frac{2!(1+|\xi|)}{(1-|\xi|^2)^2},\cdots,
    \frac{(m-k+1)!(1+|\xi|)^{m-k}}{(1-|\xi|^2)^{m-k+1}}\right)\left(2\log\frac{1}{|f(\xi)|}\right)^k.
	\end{split}\end{equation}
A direct calculation shows
$$\frac{(m-k+1)!(1+|\xi|)^{m-k}}{(1-|\xi|^2)^{m-k+1}}=\frac{(m-k+1)!}{(1-|\xi|)^{m-k+1}(1+|\xi|)}.$$
Therefore, we can deduce from (\ref{pan5}), (\ref{pan3}) and (\ref{pan4}) that
 \begin{align*}
	\left|f^{(m)}(\xi)\right|&\leq \frac{|f(\xi)|}{\left(1-\left|\xi\right|\right)^{m}}
\sum\limits_{k=1}^{m}B_{m,k}\left(1!,2!,\cdots,(m-k+1)!\right)\left(\frac{2\log\frac{1}{|f(\xi)|}}{1+\left|\xi\right|}\right)^k\\
	&= \frac{|f(\xi)|}{\left(1-\left|\xi\right|\right)^{m}}
\sum\limits_{k=1}^{m}\frac{m!}{k!} \! \begin{pmatrix}
m-1 \\
k-1
\end{pmatrix}\left(\frac{2\log\frac{1}{|f(\xi)|}}{1+\left|\xi\right|}\right)^k
	\end{align*}
for any given $\xi\in \mathbb{D}$. The proof is complete.

 \begin{remark}\label{lemma2}
    When $m=1$, the above inequality states that any holomorphic function $f$ from the unit disk into the punctured unit disk is a contraction with respect to their corresponding Poincar$\acute{\text{e}}$ metrics. That is
    $$
    \frac{|f'(z)|}{2|f(z)|\log\frac{1}{|f(z)|}}\leq\frac{1}{1-|z|^2}, \ \ z\in\mathbb{D}.
    $$
     \end{remark}

     \begin{remark}\label{lemma3}
    In particular, the equality in \rm{(\ref{20251224})} holds  due to conformal invariance when $m=1$. For example, we can take $$f(z)=e^{\frac{z-1}{z+1}}.$$
    A direct calculation implies
    $$|f'(z)|=\frac{2}{|z+1|^2}e^{\frac{|z|^2-1}{|z+1|^2}}=\frac{2|f(z)|}{1-|z|^2}\log \frac{1}{|f(z)|}.$$
    \end{remark}
    Similarly,  we can also provide the following Schwarz-Pick estimate of $f$ from $\mathbb{D}$ into the unbounded domain $\mathbb{C}\setminus\overline{\D}$.
    \begin{corollary}\label{wst1}
	Let $f:\mathbb{D}\rightarrow\mathbb{C}\setminus\overline{\D}$ be a holomorhic function. Then for each $z\in\D$ and $m\geq1$,
	\begin{align*}
	\left|f^{(m)}(z)\right|\leq \frac{|f(z)|}{\left(1-\left|z\right|\right)^{m}}
\sum\limits_{k=1}^{m}\frac{m!}{k!} \! \begin{pmatrix}
m-1 \\
k-1
\end{pmatrix}\left(\frac{2\log|f(z)|}{1+\left|z\right|}\right)^k.
	\end{align*}
    \end{corollary}
    \begin{proof}
 Take $g(z)=\log f(z).$ Then similar to the proof of Theorem \ref{cw2025}, we obtain
   \begin{align*}
	\left|f^{(m)}(z)\right|\leq \frac{|f(z)|}{\left(1-\left|z\right|\right)^{m}}
\sum\limits_{k=1}^{m}\frac{m!}{k!} \! \begin{pmatrix}
m-1 \\
k-1
\end{pmatrix}\left(\frac{2\log|f(z)|}{1+\left|z\right|}\right)^k.
	\end{align*}
The proof is complete.
    \end{proof}
   Next, for any $a\in \mathbb{D}\setminus\{0\}$, we establish a version of Schwarz-Pick estimate, which extends Theorem  \ref{cw2025}.
   \begin{theorem}	Let $f: \mathbb{D}\rightarrow \mathbb{D} \setminus \{a\}$ be a holomorphic function with $a \in \mathbb{D}$. Then for any $z \in \mathbb{D}$ and $m \geq 1$,
	 \begin{align*}
	\left|f^{(m)}(z)\right| \leq \frac{(1 + |a|)}{\left(1-\left|z\right|\right)^m}\sum_{k=1}^m \frac{k!  |a|^{k-1} }{(1 - |a|)^{2k}} \left|a - f(z)\right|^k  B_{m,k}\left(\beta_1, \cdots, \beta_{m-k+1}\right),
	\end{align*}
	where
	$$
	\beta_j(z)= \sum\limits_{k=1}^{j}\frac{j!}{k!} \! \begin{pmatrix}
j-1 \\
k-1
\end{pmatrix}\left(\frac{2\log\left|\frac{1-\bar{a}f(z)}{a-f(z)}\right|}{1+\left|z\right|}\right)^k,\ \ 1\leq j\leq m-k+1.
	$$
\end{theorem}
\begin{proof}
    For given $a\in\mathbb{D}$, let $\psi_a:\mathbb{D}\longrightarrow\mathbb{D}$ be the M$\ddot{\text{o}}$bius transformation, defined by
	$$
    \psi_a(z) = \frac{a-z}{1-\overline{a}z}
    $$
	with its inverse $\psi_a^{-1}=\psi_a.$
	Take
	$$\varphi(z) = \psi_a(f(z)) = \frac{a-f(z)}{1-\overline{a}f(z)}.$$
	Then $\varphi: \D \to \D \setminus \{0\}$ is holomorphic. It follows from Theorem \ref{cw2025} that
 \begin{align*}
	\left|\varphi^{(j)}(z)\right|\leq \frac{|\varphi(z)|}{\left(1-\left|z\right|\right)^{j}}
\sum\limits_{k=1}^{j}\frac{j!}{k!} \! \begin{pmatrix}
j-1 \\
k-1
\end{pmatrix}\left(\frac{2\log\frac{1}{|\varphi(z)|}}{1+\left|z\right|}\right)^j,\ \ j\geq1.
	\end{align*}
	Set
	$$
	\beta_j(z)= \sum\limits_{k=1}^{j}\frac{j!}{k!} \! \begin{pmatrix}
j-1 \\
k-1
\end{pmatrix}\left(\frac{2\log\left|\frac{1-\bar{a}f(z)}{a-f(z)}\right|}{1+\left|z\right|}\right)^k.
	$$
	Then
	\begin{equation}\label{equ3}
	|\varphi^{(j)}(z)| \leq  \frac{|\varphi(z)|}{\left(1-\left|z\right|\right)^{j}}  \beta_j(z).
	\end{equation}
	Notice that
    $$
    f(z)=\psi_{a}^{-1}(\varphi(z))=\psi^{-1}_{a}\left(\frac{a-f(z)}{1-\overline{a}f(z)}\right).
    $$
    Thus, by Lemma \ref{2.1}, we have
	\begin{align*}\label{equ1}
	f^{(m)}(z) = \sum_{k=1}^m \left(\psi_a^{-1}\right)^{(k)}\left(\varphi(z)\right) B_{m,k}\left(\varphi'(z),\cdots,\varphi^{(m-k+1)}(z)\right).
	\end{align*}
    This, together with
	\begin{align*}
    (\psi_a^{-1})^{(k)}(w) = \frac{k! \bar{a}^{k-1} \left(|a|^{2} - 1\right)}{(1 - \bar{a} w)^{k+1}},
    \end{align*}
	implies
	\begin{equation}\label{eq2025}
	f^{(m)}(z) = \sum_{k=1}^m \frac{k!  \overline{a}^{k-1} \left(\left|a\right|^2-1 \right)}{\left(1 - \overline{a}\varphi(z)\right)^{k+1}}  B_{m,k}\left(\varphi'(z), \varphi''(z), \cdots, \varphi^{(m-k+1)}(z)\right).
	\end{equation}
	Hence, (\ref{equ3}) and (\ref{eq2025}) yield
	$$
	\left|f^{(m)}(z)\right| \leq \sum_{k=1}^m \frac{k!  \left|a\right|^{k-1} \left(1 - \left|a\right|^2\right)}{\left(1 - \left|a\right|\right)^{k+1}}  B_{m,k}\left(\frac{|\varphi(z)|}{\left(1-\left|z\right|\right)}  \beta_1(z),\cdots,\frac{|\varphi(z)|}{\left(1-\left|z\right|\right)^{m-k+1}}  \beta_{m-k+1}(z)\right).
	$$
	From (\ref{pan3}), we get
	$$
	B_{m,k}\left(\frac{|\varphi(z)|}{\left(1-\left|z\right|\right)}  \beta_1(z),\cdots,\frac{|\varphi(z)|}{\left(1-\left|z\right|\right)^{m-k+1}}  \beta_{m-k+1}(z)\right) = \frac{\left|\varphi(z)\right|^k }{\left(1-\left|z\right|\right)^m} B_{m,k}\left(\beta_1, \cdots, \beta_{m-k+1}\right).
	$$
	Thus,
	\begin{equation}\label{equ4}
	|f^{(m)}(z)| \leq  \frac{(1 - |a|^2)}{\left(1-\left|z\right|\right)^m}\sum_{k=1}^m \frac{k!  |a|^{k-1} }{(1 - |a|)^{k+1}}  |\varphi(z)|^k  B_{m,k}(\beta_1, \cdots, \beta_{m-k+1}).
	\end{equation}
	Therefore, by (\ref{equ4}) and the fact that
	$$
	\left|\varphi(z)\right| = \left|\frac{a - f(z)}{1 - \overline{a}f(z)}\right| \leq \frac{\left|a - f(z)\right|}{1 - \left|a\right|},
	$$
	we obtain
    \begin{align*}
	\left|f^{(m)}(z)\right| &\leq  \frac{(1 - |a|^2)}{\left(1-\left|z\right|\right)^m}\sum_{k=1}^m \frac{k!  |a|^{k-1} }{(1 - |a|)^{k+1}} \frac{|a - f(z)|^k}{(1 - |a|)^k}  B_{m,k}\left(\beta_1, \cdots, \beta_{m-k+1}\right)\\
    &\leq  \frac{(1 - |a|^2)}{\left(1-\left|z\right|\right)^m}\sum_{k=1}^m \frac{k!  |a|^{k-1}   \left|a - f(z)\right|^k}{(1 - |a|)^{2k + 1}}  B_{m,k}\left(\beta_1, \cdots, \beta_{m-k+1}\right)\\
    &= \frac{ (1 + |a|)}{\left(1-\left|z\right|\right)^m}\sum_{k=1}^m \frac{k!  |a|^{k-1}}{(1 - |a|)^{2k}} \left|a - f(z)\right|^k  B_{m,k}\left(\beta_1, \cdots, \beta_{m-k+1}\right).
	\end{align*}
 The proof is complete.
    \end{proof}
	\begin{remark}\label{2026}
	When $m=1$, we have
	$$|f'(z)| \leq \frac{2(1 + |a|)}{(1 - |a|)^{2}(1 - |z|^{2})} |a - f(z)| \log \left| \frac{1 - \bar{a} f(z)}{a - f(z)} \right|, \ \ z\in\mathbb{D}.$$
	 In addition, if $a=0$, then
 \begin{align*}
	\left|f^{(m)}(z)\right| \leq \frac{1}{\left(1-\left|z\right|\right)^m} \left|f(z)\right|  B_{m,1}\left(\beta_1, \cdots, \beta_{m}\right),
	\end{align*}
	where
	$$
	\beta_j(z)= \sum\limits_{k=1}^{j}\frac{j!}{k!} \! \begin{pmatrix}
j-1 \\
k-1
\end{pmatrix}\left(\frac{2\log\frac{1}{\left|f(z)\right|}}{1+\left|z\right|}\right)^k,\ \ 1\leq j\leq m.
	$$	
By the definition of Bell polynomial, we have $B_{m,1}\left(\beta_1, \cdots, \beta_{m}\right)=\beta_m.$ Hence, it holds that
 \begin{align*}
	\left|f^{(m)}(z)\right|\leq \frac{|f(z)|}{\left(1-\left|z\right|\right)^{m}}
\sum\limits_{k=1}^{m}\frac{m!}{k!} \! \begin{pmatrix}
m-1 \\
k-1
\end{pmatrix}\left(\frac{2\log\frac{1}{|f(z)|}}{1+\left|z\right|}\right)^k,
	\end{align*}
    which reduces to Theorem \ref{cw2025}.
\end{remark}
Finally, we study the case of finitely many points in $\mathbb{D}$, which also generalizes Theorem \ref{cw2025}. For simplicity, we only need to prove the case of $f'(z)$.
\begin{corollary}
	Let $f:\mathbb{D} \rightarrow \mathbb{D}\setminus \{a_1,\cdots,a_m\}$ be a holomorphic function with $a_j\in\mathbb{D}$ for $1\leq j\leq m$. Then
	$$|f'(z)| \leq \frac{2}{1 - |z|^{2}} \left[\prod_{j=1}^{m}\frac{(1 + |a_j|)|a_j - f(z)|}{(1 - |a_j|)^{2}}  \log \left| \frac{1 - \overline{a_j} f(z)}{a_j - f(z)} \right|\right]^{1/ m}, \ \ z\in\mathbb{D}.$$
\end{corollary}
\begin{proof}
From Remark \ref{2026}, for each $a_j\in\mathbb{D}$, we have
$$|f'(z)| \leq \frac{2(1 + |a_j|)}{(1 - |a_j|)^{2}(1 - |z|^{2})} |a_j - f(z)| \log \left| \frac{1 - \overline{a_j} f(z)}{a_j - f(z)} \right|, \ \ z\in\mathbb{D},$$
for $j=1,\cdots,m$. This implies
$$|f'(z)| \leq \frac{2}{1 - |z|^{2}} \left[\prod_{j=1}^{m}\frac{(1 + |a_j|)|a_j - f(z)|}{(1 - |a_j|)^{2}}  \log \left| \frac{1 - \overline{a_j} f(z)}{a_j - f(z)} \right|\right]^{1/m}.$$
The proof is complete.
\end{proof}
With the help of the geometric mean inequality, we also give another form of the Schwarz-Pick estimates of holomorphic functions that avoid any finitely many points of the unit disk.
\begin{corollary}
	Let $f:\mathbb{D} \rightarrow \mathbb{D}\setminus \{a_1,\cdots,a_m\}$ be a holomorphic function with $a_j\in\mathbb{D}$ for $1\leq j\leq m$. Then
    $$\left|f'(z)\right| \leq \frac{2}{m\left(1 - |z|^{2}\right)}\sum\limits_{j=1}^{m}\frac{\left(1 + |a_j|\right)\left|a_j - f(z)\right|}{\left(1 - |a_j|\right)^{2}}\log \left| \frac{1 - \overline{a_j} f(z)}{a_j - f(z)} \right|,\ \  z\in\mathbb{D}.$$
\end{corollary}
\begin{proof}
For each $j=1,2,\cdots,m$ and $z\in\mathbb{D}$, it follows from Remark \ref{2026} that
$$|f'(z)| \leq \frac{2(1 + |a_j|)}{(1 - |a_j|)^{2}(1 - |z|^{2})} |a_j - f(z)| \log \left| \frac{1 - \overline{a_j} f(z)}{a_j - f(z)} \right|.$$
Then the geometric mean inequality
$$
\sqrt[m]{x_1 \cdots x_m}\leq \frac{x_1+\cdots+x_m}{m},
$$
shows that for any $z\in\mathbb{D}$,
$$\left|f'(z)\right| \leq \frac{2}{m\left(1 - |z|^{2}\right)}\sum\limits_{j=1}^{m}\frac{\left(1 + |a_j|\right)\left|a_j - f(z)\right|}{\left(1 - |a_j|\right)^{2}}\log \left| \frac{1 - \overline{a_j} f(z)}{a_j - f(z)} \right|.$$
The proof is complete.
\end{proof}

\section{Schwarz-Pick Estimates on the unit ball}
In this section, we establish a new Schwarz-Pick estimate for nonvanishing holomorphic functions on $\mathbb{B}_n$, which is a generalization of Theorem \ref{cw2025}.
Let
$$\mathbb{D}_{z_0,r}=\{z\in\mathbb{C}:|z-z_0|<r\}$$
and
$$\mathbb{D}^n=\{z=(z_1,\cdots,z_n)\in\mathbb{C}^n:|z_j|<1, 1\leq j\leq n\}$$
be a disk with center at $z_0$ and radius $r$ in $\mathbb{C}$ and the unit polydisk in $\mathbb{C}^n$, respectively.
\begin{lemma}[\cite{DCP2010}]\label{pjq}
	For two given points $p,q\in\mathbb{B}_n$ with $q\neq p$, let $L(z)=p+z(q-p)$ for $z\in\mathbb{C}$. Then $L(\mathbb{D}_{c_{p,q},r_{p,q}})\subset \mathbb{B}_n$, where
\begin{center}
$c_{p,q} = -\frac{\langle p,q-p\rangle}{|q-p|^2}$,\  $r_{p,q}=\sqrt{\frac{1-|p|^2}{|q-p|^2}+\left|\frac{\langle p,q-p\rangle}{|q-p|^2}\right|^2}.$
\end{center}
\end{lemma}
\begin{lemma}\label{pjq1}
Let $g:\mathbb{D}_{z_0,r}\rightarrow\mathbb{D}\setminus\left\{0\right\}$ be a holomorphic function. Then for any $z\in D_{z_0,r}$ and $m\geq 1$,
 \begin{align*}
	\left|g^{(m)}(z)\right|\leq \frac{|g(z)|}{\left(r-\left|z-z_{0}\right|\right)^{m}}
\sum\limits_{k=1}^{m}\frac{m!}{k!} \! \begin{pmatrix}
m-1 \\
k-1
\end{pmatrix}\left(\frac{2r\log\frac{1}{|g(z)|}}{r+\left|z-z_{0}\right|}\right)^k.
\end{align*}
\end{lemma}
\begin{proof}
	Take
    $$
    \varphi(z)=g(rz+z_0),\ \ z\in\mathbb{D}.
    $$
 Then $\varphi:\mathbb{D}\rightarrow\mathbb{D}\setminus\{0\}$ is a holomorphic function. It follows from Theorem \ref{cw2025} that
 \begin{align*}
	\left|\varphi^{(m)}(z)\right|\leq \frac{|\varphi(z)|}{\left(1-\left|z\right|\right)^{m}}
\sum\limits_{k=1}^{m}\frac{m!}{k!} \! \begin{pmatrix}
m-1 \\
k-1
\end{pmatrix}\left(\frac{2\log\frac{1}{|\varphi(z)|}}{1+\left|z\right|}\right)^k
	\end{align*}
	for any $z\in\mathbb{D}$ and $m\geq1$. Also, note that $\varphi^{(m)}(\frac{\zeta-z_{0}}{r})=r^{m}g^{(m)}(\zeta)$ for any $\zeta \in D_{z_0,r}$.
	Setting $z=\frac{\zeta-z_{0}}{r}$ gives
\begin{align*}
	\left|r^{m}g^{(m)}(\zeta)\right|\leq \frac{|g(\zeta)|}{\left(1-\left|\frac{\zeta-z_{0}}{r}\right|\right)^{m}}
\sum\limits_{k=1}^{m}\frac{m!}{k!} \! \begin{pmatrix}
m-1 \\
k-1
\end{pmatrix}\left(\frac{2\log\frac{1}{|g(\zeta)|}}{1+\left|\frac{\zeta-z_{0}}{r}\right|}\right)^k.
	\end{align*}
By a direct computation, we get
	\begin{align*}
	\left|g^{(m)}(\zeta)\right|\leq \frac{|g(\zeta)|}{\left(r-\left|\zeta-z_{0}\right|\right)^{m}}
\sum\limits_{k=1}^{m}\frac{m!}{k!} \! \begin{pmatrix}
m-1 \\
k-1
\end{pmatrix}\left(\frac{2r\log\frac{1}{|g(\zeta)|}}{r+\left|\zeta-z_{0}\right|}\right)^k.
	\end{align*}
The proof is complete.
\end{proof}
\textbf{Proof of Theorem \ref{pjq2026-3}.} For given $p,q\in \mathbb{B}_n$ with $q\neq p$, take
    $$g\left(z\right)=f\left(L(z)\right)
    $$
    for $z\in\mathbb{C}$, where $L(z)=p+z(q-p)$. By Lemma \ref{pjq}, $g:\mathbb{D}^{*}\longrightarrow\mathbb{D}\setminus\left\{0\right\}$ is a holomorphic function, where
    \begin{center}
$\mathbb{D}^{*}=\left\{z\in\C:\left|z-c_{p,q}\right|<r_{p,q}\right\}$,\ $c_{p,q} = -\frac{\langle p,q-p\rangle}{|q-p|^2}$,\  $r_{p,q}=\sqrt{\frac{1-|p|^2}{|q-p|^2}+\left|\frac{\langle p,q-p\rangle}{|q-p|^2}\right|^2}.$
    \end{center}
	Applying Lemma \ref{pjq1} at $z=0\in\mathbb{D}^*$, we derive
\begin{equation}\label{cw2026}\begin{split}
	\left|g^{(m)}(0)\right|\leq \frac{|g(0)|}{\left(r_{p,q}-\left|c_{p,q}\right|\right)^{m}}
\sum\limits_{k=1}^{m}\frac{m!}{k!} \! \begin{pmatrix}
m-1 \\
k-1
\end{pmatrix}\left(\frac{2r_{p,q}\log\frac{1}{|g(0)|}}{r_{p,q}+\left|c_{p,q}\right|}\right)^k
	 \end{split}\end{equation}
	for any $m\geq1$. Note that
	$$g(0)=f(p),\ \ g^{(m)}(0)=\sum\limits_{|v|=m}\frac{m!}{v!}\frac{\partial^m f(p)}{\partial w_{1}^{v_{1}}\cdots\partial w_{n}^{v_{n}}}(q-p)^{v}=D^mf(p)\cdot(q-p)^m.$$
    Here, $v=(v_1, \cdots, v_n)$ is a multi-index. This, together with (\ref{cw2026}), shows
    \begin{equation}\label{cw2027}\begin{split}
	\left|D^mf(p)\cdot\beta^m\right|\leq \frac{|f(p)|}{\left(\sqrt{A}-\left|\langle p,\beta\rangle\right|\right)^{m}}
    \sum\limits_{k=1}^{m}\frac{m!}{k!} \! \begin{pmatrix}
    m-1 \\
    k-1
    \end{pmatrix}\left(\frac{2\sqrt{A}\log\frac{1}{|f(p)|}}{\sqrt{A}+\left|\langle p,\beta\rangle\right|}\right)^k,
	 \end{split}\end{equation}
    where $A=1-\left|p\right|^2+\left|\langle p,\beta\rangle\right|^2$ and $\beta=\frac{q-p}{\left|q-p\right|}$. A direct computation yields
    \begin{align*}
    \frac{\left(\sqrt{A}\right)^k}{\left(\sqrt{A}-\left|\langle p,\beta\rangle\right|\right)^m\left(\sqrt{A}+\left|\langle p,\beta\rangle\right|\right)^k}&=\frac{1}{\left(\sqrt{A}-\left|\langle p,\beta\rangle\right|\right)^{m-k}}\left(\frac{\sqrt{A}}{1-\left|p\right|^2}\right)^k\\
    &=\frac{\left(\sqrt{A}+\left|\langle p,\beta\rangle\right|\right)^{m-k}}{\left(1-|p|^2\right)^{m-k}}H_{p}(\beta,\beta)^{k/2},
    \end{align*}
     where
    $$ H_p(\beta,\beta)= \frac{(1-|p|^2)|\beta|^2 + |\langle\beta,p\rangle|^2}{(1-|p|^2)^2}$$
    is the Bergman metric in the direction of $\beta\in \partial \mathbb{B}_n$.

    Notice that
     $$ H_p(\beta,\beta)= \frac{(1-|p|^2)|\beta|^2 + |\langle\beta,p\rangle|^2}{(1-|p|^2)^2}\leq \frac{1}{(1-|p|^2)^2}.$$
    Hence,
      \begin{equation}\label{cw2035}
    \frac{\left(\sqrt{A}\right)^k}{\left(\sqrt{A}-\left|\langle p,\beta\rangle\right|\right)^m\left(\sqrt{A}+\left|\langle p,\beta\rangle\right|\right)^k}\leq \frac{\left(\sqrt{A}+\left|\langle p,\beta\rangle\right|\right)^{m-k}}{\left(1-|p|^2\right)^{m}}.
    \end{equation}
   Therefore, (\ref{cw2027}) and (\ref{cw2035}) show
	 \begin{align*}\label{cw2029}
	\left|D^mf(p)\cdot\beta^m\right|\leq |f(p)|\left(\frac{\sqrt{A}+\left|\langle p,\beta\rangle\right|}{1-|p|^2}\right)^m
    \sum\limits_{k=1}^{m}\frac{m!}{k!} \! \begin{pmatrix}
    m-1 \\
    k-1
    \end{pmatrix}\left(\frac{2\log\frac{1}{|f(p)|}}{\sqrt{A}+\left|\langle p,\beta\rangle\right|}\right)^k.
	 \end{align*}
Since $\beta$ is an arbitrary unit vector, this means that (\ref{cw2032}) holds for any $p\in\mathbb{B}_n$ and $\beta\in\partial\mathbb{B}_n$.
The proof is complete.

 \begin{remark} When $z=\mathbf{0}\in\mathbb{B}_n$, it holds that for any $\beta\in\partial\mathbb{B}_n$,
	\begin{align*}
		\left|D^mf(\mathbf{0})\cdot\beta^m\right|\leq |f(\mathbf{0})|
    \sum\limits_{k=1}^{m}\frac{m!}{k!} \! \begin{pmatrix}
    m-1 \\
    k-1
    \end{pmatrix}\left(2\log\frac{1}{|f(\mathbf{0})|}\right)^k,\ \ m\geq1.
	\end{align*}
If $n=1$,  then Theorem {\rm\ref{pjq2026-3}} reduces to Theorem {\rm\ref{cw2025}}.
\end{remark}
Now we apply Theorem \ref{pjq2026-3} to obtain the following interesting examples which tries to find higher order Schwarz-Pick type estimates of holomorphic mappings from $\mathbb{B}_n$ into the multiply connected domain in $\mathbb{C}^n$.
\begin{example}\label{pjq2026-18}
	 Let $\varphi=(\varphi_1,\varphi_2,\cdots,\varphi_n):\mathbb{B}_n\rightarrow\mathbb{D}^n\setminus\left\{\prod\limits_{j=1}^{n}z_j=0\right\}$ be a holomorphic mapping.
 Take $\Phi(z)=\prod\limits_{j=1}^{n}\varphi_j(z)$. Then $\Phi(z): \mathbb{B}_n\rightarrow\mathbb{D}\setminus\{0\}$ is a holomorphic function. By Theorem {\rm\ref{pjq2026-3}},
  for any $z\in\mathbb{B}_n$ and $\beta\in\partial \mathbb{B}_n$,
		\begin{align*}
		\left|D^m\Phi(z)\cdot\beta^m\right|
\leq |\Phi(z)|\left(\frac{\sqrt{A}+\left|\langle z,\beta\rangle\right|}{1-|z|^2}\right)^m
    \sum\limits_{k=1}^{m}\frac{m!}{k!} \! \begin{pmatrix}
    m-1 \\
    k-1
    \end{pmatrix}\left(\frac{2\sum\limits_{j=1}^n\log\frac{1}{|\varphi_j(z)|}}{\sqrt{A}+\left|\langle z,\beta\rangle\right|}\right)^k,
		\end{align*}
where $A=1-\left|z\right|^2+\left|\langle z,\beta\rangle\right|^2$ and $m\geq1$.
\end{example}

\begin{example}\label{pjq2026-19}
	 Let $\varphi=(\varphi_1,\varphi_2,\cdots,\varphi_n):\mathbb{B}_n\rightarrow\mathbb{B}_n\setminus\left\{\sum\limits_{j=1}^{n}z^2_j=0\right\}$ be a holomorphic mapping.  Take $\Phi(z)=\sum\limits_{j=1}^{n}\varphi^2_j(z)$. Then $\Phi(z): \mathbb{B}_n\rightarrow\mathbb{D}\setminus\{0\}$ is a holomorphic function.  By Theorem {\rm\ref{pjq2026-3}},
  for any $z\in\mathbb{B}_n$ and $\beta\in\partial \mathbb{B}_n$,
		\begin{align*}
		\left|D^m\Phi(z)\cdot\beta^m\right|
\leq |\Phi(z)|\left(\frac{\sqrt{A}+\left|\langle z,\beta\rangle\right|}{1-|z|^2}\right)^m
    \sum\limits_{k=1}^{m}\frac{m!}{k!} \! \begin{pmatrix}
    m-1 \\
    k-1
    \end{pmatrix}\left(\frac{2\log\frac{1}{\left|\Phi(z)\right|}}{\sqrt{A}+\left|\langle z,\beta\rangle\right|}\right)^k,
		\end{align*}
where $A=1-\left|z\right|^2+\left|\langle z,\beta\rangle\right|^2$ and $m\geq1$.
\end{example}

 Unfortunately, since the multiply connected domain in $\mathbb{C}^n$ is completely different from that of $\mathbb{C}$,
 we do not know whether similar results hold or not. We list it as an open problem.
\begin{problem}
Let $\Omega\subset\mathbb{C}^n$ be a multiply connected bounded domain. For a holomorphic mapping $f:\mathbb{B}_n\rightarrow\Omega$, is it possible to establish higher order Schwarz-Pick type estimates analogous to those presented in Theorem {\rm\ref{cw2025}?}
\end{problem}

\section{Application}
In this section, we apply Schwarz-Pick estimates of nonvanishing holomorphic functions to obtain the $\L$ojasiewicz inequality in the case of several complex variables, which describes the connection between holomorphic functions and their zero sets. To proceed further, we first introduce some definitions and lemmas.
  \begin{definition}
    If $V_1$,$\cdots$,$V_{k}$ and $W$ are linear spaces, a mapping
    $$A: V_1\times V_2\times \cdots\times V_k\rightarrow W$$
    is called $k$-multilinear if $A(e_1,\cdots,e_k)$ is linear in each argument separately. Linearity in the $j$th argument means that
    $$A(e_1,\cdots,\alpha e_j+\beta \widehat{e}_j,\cdots,e_k)=\alpha A(e_1,\cdots,e_j,\cdots,e_k)+\beta A(e_1,\cdots,\widehat{e}_j,\cdots,e_k)$$
    for any $e_j,\widehat{e}_j\in V_j\ (1\leq j\leq k)$ and all scalars $\alpha$, $\beta$.
    \end{definition}
\begin{lemma}[\cite{AM88}]\label{symmetric}
    Let $V$ and $W$ be two Banach spaces with norms $\|\cdot\|_V$ and $\|\cdot\|_{W}$. Suppose that
 $$
A: \underbrace{V \times V \times \cdots \times V}_{k \text{ copies of } V} \rightarrow W
$$
    is a symmetric and $k$-multilinear mapping. The symmetry means that for any permutation $\sigma$ of $\{1, \cdots, k\}$,
    $$ A(v_{\sigma(1)}, \cdots, v_{\sigma(k)}) = A(v_1, \cdots, v_k).$$
    Then
    \begin{equation} \label{eq:main_inequality}
    \sup_{\|e_1\|_V=1, \cdots, \|e_k\|_V=1} \|A(e_1, \cdots, e_k)\|_W \leq \frac{k^k}{k!} \sup_{\|e\|_V= 1} \|A(e, \cdots, e)\|_W.
    \end{equation}
\end{lemma}
\begin{lemma}\label{cw12}
For any $\varepsilon>0$ and $0<x<1$, we have $\log \frac{1}{x}\leq \frac{1}{\varepsilon} x^{-\varepsilon}$.
\end{lemma}
\begin{proof}
Note that $t^{-1}\leq t^{-\varepsilon-1}$ for any $0<t<1$. Then
$$
\log \frac{1}{x}=\int_{x}^{1}t^{-1}dt\leq \int_{x}^{1}t^{-1-\varepsilon}dt=-\frac{1}{\varepsilon}t^{-\varepsilon}\bigg|_{t=x}^{t=1}=\frac{x^{-\varepsilon}-1}{\varepsilon}\leq \frac{1}{\varepsilon} x^{-\varepsilon}.$$
The proof is complete.
\end{proof}
Applying Theorem \ref{cw2025} at $z=0$, we have the following lemma.
  \begin{lemma}\label{lemma2}
	Let $f:\mathbb{D}\rightarrow\mathbb{D}\setminus\left\{0\right\}$ be a holomorphic function. Then for any $m\geq 1$,
	\begin{equation}\label{cw}
	\left|f^{(m)}(0)\right|\leq\left|f(0)\right|\sum_{k=1}^{m}\frac{m!}{k!} \! \begin{pmatrix}
m-1 \\
k-1
\end{pmatrix}\left(2\log\frac{1}{|f(0)|}\right)^{k}.
	\end{equation}
    \end{lemma}


	\textbf{Proof of Theorem \ref{wst2030}.} Without loss of generality, we may assume that $q=\mathbf{0}$. Then $\mathbf{0}\in Z_f$. Since $f$ has vanishing order $k$
 at $\mathbf{0}$, we can choose $1\leq i_1,\cdots,i_k\leq n$ such that
	$$\frac{\partial^k f}{\partial z_{i_1}\cdots\partial z_{i_k}}(\mathbf{0})\neq 0.$$
	Then there exists some neighborhood $U_{\mathbf{0}}$ of $\mathbf{0}$ and a constant $c_{0}>0$ such that
	\begin{equation}\label{cw34}
    \left|\frac{\partial^k f}{\partial z_{i_1}\cdots\partial z_{i_k}}(z)\right|\geq c_0,\quad z\in U_\mathbf{0}.
    \end{equation}For small enough $r>0$, we have $|f(z)|<1$ on $B(\mathbf{0}, r)$, where $B(\mathbf{0}, r)=\{z\in\mathbb{C}^n:|z|<r\}\subset U_{\mathbf{0}}$.  We may assume that there exists $p\in B\left(\mathbf{0},\frac{r}{2}\right)$ such that $p\notin Z_f$. Let $z_{p}\in Z_f$ be the point satisfying
    $$
    \operatorname{dist}(p, Z_f)=\inf_{z\in Z_f}\left|p-z\right|=\left|p-z_{p}\right|.
    $$
    Set $d_{p}=\operatorname{dist}(p, Z_f)$. By the definition of $\operatorname{dist}(p, Z_f)$ and $\mathbf{0}\in Z_f$, we have $d_{p}\leq|p|$. For any $\omega=(\omega_1,\cdots,\omega_n)\in \C^{n}$ with $|\omega|=1$, we define a complex line
	$$\gamma_\omega(\zeta)=p+\left(d_p \omega\right)\zeta$$
    on  $\mathbb{D}$.	Then
    \begin{center}
    $\gamma_{\omega}(0)=p$ and $\gamma_{\omega}^{\prime}(0)=d_{p} \omega$.
    \end{center}
    Observe that for any $\zeta\in\mathbb{D}$,
	\begin{align*}
	\left|\gamma_\omega(\zeta)\right|&\leq|p|+d_p|\omega||\zeta|< 2|p|< 2\cdot\frac{r}{2}=r.
	\end{align*}
	Take $g(\zeta)=f\left(\gamma_{\omega}(\zeta)\right)$.  Since $|f(z)|<1$ for any $z\in B(\mathbf{0}, r)$, we know that $g: \D\rightarrow \D$ is holomorphic. Notice that
	$$\left|\gamma_\omega(\zeta)-p\right|=\left|\left(d_p \omega\right)\zeta\right|<d_p$$
    for any $\zeta\in\mathbb{D}$. Then this implies $g(\zeta)\neq 0$ for any $\zeta\in\mathbb{D}$.
	A direct computation shows
	\begin{equation}\label{equation1}
    g^{\prime}(\zeta)=\sum_{j=1}^n\frac{\partial f}{\partial z_j}\left(\gamma_\omega(\zeta)\right) \omega_j d_p
    \end{equation}
	with $\gamma_{\omega}(\zeta)=(\gamma_{\omega,1}(\zeta),\cdots,\gamma_{\omega,n}(\zeta))$.
    By Lemma \ref{lemma2}, we get
	$$\left|g^{\prime}(0)\right|\leq 2|g(0)|\log\frac{1}{|g(0)|}.$$
     This, together with (\ref{equation1}) and $g(0)=f(\gamma_{\omega}(0))=f(p)$, gives
	$$|g'(0)|=\left|\sum_{j=1}^{n}\frac{\partial f}{\partial z_{j}}(p) \omega_{j}\right| d_{p}\leq 2|f(p)|\log\frac{1}{|f(p)|}.$$
	Choosing $\omega$ so that $\omega_{j}=1$ for some $j$ and $\omega_{k}=0$ for any $k\neq j$, we obtain
	$$\left|\frac{\partial f}{\partial z_j}(p)\right| d_p\leq 2|f(p)|\log\frac{1}{|f(p)|}$$
	for $j=1,\cdots,n$.
	Note that
	\begin{equation}\label{equation2} g^{\prime\prime}(\zeta)=\sum_{i, j=1}^n\frac{\partial^2 f}{\partial z_i\partial z_j}\left(\gamma_\omega(\zeta)\right)\frac{\partial\gamma_{\omega,i}}{\partial\zeta}\frac{\partial\gamma_{\omega,j}}{\partial\zeta}=\sum_{i, j=1}^n\frac{\partial^2 f}{\partial z_i\partial z_j}\left(\gamma_\omega(\zeta)\right) \omega_i \omega_j d_p^2.
    \end{equation}
	It follows from Lemma \ref{lemma2} that
	\begin{equation}\label{equation3}
    \left|g^{\prime\prime}(0)\right|\leq 4|g(0)|\left(\log\frac{1}{|g(0)|}+\left(\log\frac{1}{|g(0)|}\right)^{2}\right).
    \end{equation}
	Combining (\ref{equation2}) and (\ref{equation3}), we have
	\begin{equation}\label{equation4}
	|g''(0)|=\left|\sum_{i, j=1}^n\frac{\partial^2 f}{\partial z_i\partial z_j}(p) \omega_i \omega_j\right| d_p^2\leq 4|f(p)|\left(\log\frac{1}{|f(p)|}+\left(\log\frac{1}{|f(p)|}\right)^2\right).
	\end{equation}
	Define the mapping
	$$(\lambda, \mu)\mapsto\sum_{i, j=1}^n\frac{\partial^2 f}{\partial z_i\partial z_j}(p) \lambda_i \mu_j.$$
	Then we know that this is a symmetric and 2-multilinear mapping. Hence, by Lemma \ref{symmetric},
	\begin{equation}\label{ptc}\sup_{\|\lambda\|=1,\|\mu\|=1}\left|\sum_{i, j=1}^n\frac{\partial^2 f}{\partial z_i\partial z_j}(p) \lambda_i \mu_j\right|\leq 2\sup_{\|\lambda\|=1}\left|\sum_{i, j=1}^n\frac{\partial^2 f}{\partial z_i\partial z_j}(p) \lambda_i \lambda_j\right|.\end{equation}
	Choosing $\lambda=(0,\cdots,\overbrace{1}^{i},\cdots, 0)$ and $\mu=(0,\cdots,\overbrace{1}^{j},\cdots, 0)$ in (\ref{ptc}) and applying (\ref{equation4}), we get
	\begin{align*}
    \left|\frac{\partial^2 f}{\partial z_i\partial z_j}(p)\right| d_p^2&\leq\sup_{\|\lambda\|=1,\|\mu\|=1}\left|\sum_{i, j=1}^n\frac{\partial^2 f}{\partial z_i\partial z_j}(p) \lambda_i \mu_j\right|d_p^2\\
    &\leq 2\sup_{\|\lambda\|=1}\left|\sum_{i, j=1}^n\frac{\partial^2 f}{\partial z_i\partial z_j}(p) \lambda_i \lambda_j\right|d_p^2\\
    &\leq 8|f(p)|\left(\log\frac{1}{|f(p)|}+\left(\log\frac{1}{|f(p)|}\right)^2\right).
    \end{align*}
	In particular, taking $k$th derivatives of $g$, we obtain
	\begin{equation}\label{pjq2026}
    g^{(k)}(\zeta)=\left(\sum_{1\leq i_1,\cdots, i_k\leq n}\frac{\partial^k f}{\partial z_{i_1}\cdots\partial z_{i_k}}\left(\gamma_\omega(\zeta)\right) \omega_{i_1}\cdots \omega_{i_k}\right) d_p^k.
    \end{equation}
	Lemma \ref{lemma2} gives
	\begin{equation}\label{pjq2026-1}
    \left|g^{(k)}(0)\right|\leq C_{k}|g(0)|\left(\sum_{j=1}^{k}\left(\log\frac{1}{|g(0)|}\right)^{j}\right),
    \end{equation}
	where $C_k>0$ is a constant with depending on $k$. Combining (\ref{pjq2026}) and (\ref{pjq2026-1}), we have
	$$|g^{(k)}(0)|=\left|\sum_{1\leq i_{1},\cdots, i_{k}\leq n}\frac{\partial^{k} f}{\partial z_{i_{1}}\cdots\partial z_{i_{k}}}(p) \omega_{i_{1}}\cdots \omega_{i_{k}}\right| d_{p}^{k}\leq C_{k}|f(p)|\left(\sum_{j=1}^{k}\left(\log\frac{1}{|f(p)|}\right)^{j}\right).$$
	The  ‌arbitrariness‌ of $\omega$ implies
	\begin{equation}\label{wst2025}
\sup_{\|\omega\|=1}\left|\sum_{1\leq i_{1},\cdots, i_{k}\leq n}\frac{\partial^{k} f}{\partial z_{i_{1}}\cdots\partial z_{i_{k}}}(p) \omega_{i_{1}}\cdots \omega_{i_{k}}\right| d_{p}^{k}\leq C_{k}|f(p)|\left(\sum_{j=1}^{k}\left(\log\frac{1}{|f(p)|}\right)^{j}\right).
   \end{equation}
	Furthermore, set the mapping
	\begin{equation}\label{pjq2027-1}
\left(u_{1},\cdots, u_{k}\right)\mapsto\sum_{1\leq i_{i},\cdots, i_{k}\leq n}\frac{\partial^{k} f}{\partial z_{i_{1}}\cdots\partial z_{i_{k}}}(p)\left(u_{1}\right)_{i_{1}}\cdots\left(u_{k}\right)_{i_{k}}.
\end{equation}
	Then we know that the mapping in (\ref{pjq2027-1}) is a symmetric and $k$-multilinear mapping on $\C^{n}$ with $u_j=((u_j)_1,(u_j)_2,\cdots,(u_j)_n)$ for $1\leq j\leq k$. Hence, by Lemma 4.2 and (\ref{wst2025}), we get
	\begin{equation}\label{wst2026}\begin{split}
	&\sup_{\left\|u_j\right\|=1,1\leq j\leq k}\left|\sum_{1\leq i_1,\cdots, i_k\leq n}\frac{\partial^k f}{\partial z_{i_1}\cdots\partial z_{i_k}}(p)\left(u_1\right)_{i_1}\cdots\left(u_k\right)_{i_k}\right| d_p^k\\
	&\quad\leq\frac{k^k}{k!}\sup_{\|\omega\|=1}\left|\sum_{1\leq i_1,\cdots, i_k\leq n}\frac{\partial^k f}{\partial z_{i_1}\cdots\partial z_{i_k}}(p) \omega_{i_1}\cdots \omega_{i_k}\right| d_p^k\\
	&\quad\leq C_k\frac{k^k}{k!}|f(p)|\left(\sum_{j=1}^k\left(\log\frac{1}{|f(p)|}\right)^j\right).
	\end{split}\end{equation}
    We choose $u_{j}=(0,\cdots,\overbrace{1}^{i_{j}},\cdots, 0)$ in (\ref{wst2026}) to get
	\begin{equation}\label{equation5}
	\left|\frac{\partial^{k} f}{\partial z_{i_{1}}\cdots\partial z_{i_{k}}}(p)\right|d_{p}^{k}\leq C_k\frac{k^k}{k!}|f(p)|\left(\sum_{j=1}^{k}\left(\log\frac{1}{|f(p)|}\right)^{j}\right).
	\end{equation}
    For any $\varepsilon\in(0,1)$, by Lemma \ref{cw12},
    $$\log \frac{1}{|f(p)|}\leq \frac{k}{\varepsilon}|f(p)|^{-\frac{\varepsilon}{k}}.$$
     Then it holds that
    \begin{equation}\label{pjq2026-2}
    \sum_{j=1}^k\left(\log\frac{1}{|f(p)|}\right)^j\leq \sum_{j=1}^k\left(\frac{k}{\varepsilon}\right)^{j}|f(p)|^{-\frac{j\varepsilon}{k}}\leq k\left(\frac{k}{\varepsilon}\right)^k|f(p)|^{-\varepsilon}.
    \end{equation}
	 Therefore, together with (\ref{cw34}), (\ref{equation5}) and (\ref{pjq2026-2}), we have
	\begin{align*}\operatorname{dist}(p, Z_f)^k&\leq \frac{C_k}{c_0}\frac{k^k}{k!}|f(p)|\left(\sum_{j=1}^k\left(\log\frac{1}{|f(p)|}\right)^j\right)\\
&\leq \left(\frac{k}{\varepsilon}\right)^k\frac{C_k}{c_0}\frac{k^k}{(k-1)!}|f(p)|^{1-\varepsilon}.
\end{align*}
     This means
    $$|f(p)|\geq C\operatorname{dist}(p, Z_f)^\frac{k}{1-\varepsilon},$$
   where $C=\left[\left(\frac{k}{\varepsilon}\right)^k\frac{C_k}{c_0}\frac{k^k}{(k-1)!}\right]^{\frac{1}{\varepsilon-1}}$.
 The proof is complete.

\begin{remark}
For the case of one complex variable, notice that the zero set consists of all isolated points. Then we obtain a better inequality
$$
|f(z)| \geq C \, \mathrm{dist}(z, Z_f)^k,
$$
where the constant $C > 0$ only depends on $f$. In fact, we may assume that
$$f(z)=(z-q)^kh(z),$$
where $h(z)$ is holomorphic on a neighborhood $U$ of $q$ and $h(z)\neq0$ on $U$. This implies
$$|f(z)|= |z-q|^k |h(z)|\geq C\, \mathrm{dist}(z, Z_f)^k$$
for any $z\in U$.
\end{remark}
In what follows, rather than comparing the relative orders of $\log \frac{1}{|f(z)|}$ and $|f(z)|$, we obtain a more direct result that includes the logarithmic term to characterize the connection between holomorphic functions and $\mathrm{dist}(z,Z_f)$.
\begin{corollary}
Let $\Omega\subset\mathbb{C}^n\ (n\geq2)$ be a domain with $q\in\Omega$. Suppose that $f$ is a holomorphic function on $\Omega$ and $f(q)=0$. If $k$ is the vanishing order of $f$ at $q$, then there exists a neighborhood $U$ of $q$ such that
	\begin{align*}
	\mathrm{dist}(z,Z_f) \leq C|f(z)|^{\frac{1}{k}}\left( \log \frac{1}{|f(z)|} \right)
	\end{align*}
	for any $z\in U$, where $C>0$ is a constant that depends on $f$.
\end{corollary}
\begin{proof}
Withou loss of generality, we may assume that $q=\mathbf{0}$. Since $f$ is holomorphic and $f(\mathbf{0})=0$, there exists a neighborhood  $U$ of $\mathbf{0}$ such that $|f(z)|\leq \frac{1}{e}$ for any $z\in U$. By the proof of Theorem \ref{wst2030}, we have
\begin{align*}
 \operatorname{dist}(z, Z_f)^k&\leq \frac{C_k}{c_0}\frac{k^k}{k!}|f(z)|\left(\sum_{j=1}^k\left(\log\frac{1}{|f(z)|}\right)^j\right)\\
 &\leq C|f(z)|\left( \log \frac{1}{|f(z)|} \right)^k
 \end{align*}
 for any $z\in U$, where $C=\frac{C_k}{c_0}\frac{k^k}{(k-1)!}$. The proof is complete.
\end{proof}

{\bf Acknowledgments:} Xiaomin Tang is supported by Zhejiang Provincial Natural Science Foundation of China (No.LZ24A010004). Wei Chen is supported by Zhejiang Provincial Natural Science Foundation of China (No.LQN25A010015).

\begin{center}

\noindent
\end{center}
\end{document}